\documentclass[11pt]{article}
\usepackage{Basic,caption}
\usepackage[authoryear,round]{natbib}

\title{Spectral Gap of Metropolis Algorithms for Non-smooth Distributions under Isoperimetry}

\author{Shuigen Liu\thanks{Department of Mathematics, National University of Singapore, Singapore 119076 (\texttt{shuigen@u.nus.edu}, \texttt{mattxin@nus.edu.sg}).} \and Xin T.~Tong\footnotemark[1]}

\date{\today}

\newcommand{\CPI}{\mathsf{C_{PI}}}
\newcommand{\CLSI}{\mathsf{C_{LSI}}}

\begin{document}

\maketitle

\begin{abstract}
Metropolis algorithms are classical tools for sampling from target distributions, with broad applications in statistics and scientific computing.
Their convergence speed is governed by the spectral gap of the associated Markov operator.
Recently, \cite{24-AAP2058} derived the first explicit bounds for the spectral gap of Random--Walk Metropolis when the target distribution is smooth and strongly log-concave.
However, existing literature rarely discusses non-smooth targets.
In this work, we derive explicit spectral gap bounds for the random-walk Metropolis and Metropolis--adjusted Langevin algorithms over a broad class of non-smooth distributions.
Moreover, combining our analysis with a recent result in \cite{goyal2025mixing}, we extend these bounds to targets satisfying a Poincar\'e or log-Sobolev inequality, beyond the strongly log-concave setting.
Our theoretical results are further supported by numerical experiments.
\end{abstract}



\section{Introduction}
\label{sec:intro}

Sampling from a target distribution is a fundamental problem in statistics, scientific computing, and machine learning \citep{pmlr-v75-dwivedi18a,JMLR:v23:20-527,chewi2023logconcave}.
For decades, Markov Chain Monte Carlo (MCMC) methods \citep{MR2742422} have been the workhorse for sampling from complex distributions, among which Metropolis algorithm \citep{10.1063/1.1699114} is one of the simplest yet most widely used.
Metropolis algorithm consists of two steps:
(1) the proposal step, where a candidate sample is generated from a proposal distribution; and (2) the acceptance-rejection step, where the proposal is accepted or rejected according to an acceptance probability.
The flexibility in choosing the proposal distribution gives rise to various Metropolis algorithms, such as the Random--Walk Metropolis (RWM) \citep{MR1389882}, the Metropolis--adjusted Langevin algorithm (MALA) \citep{MR1440273}, and the Hamiltonian Monte Carlo (HMC) \citep{DUANE1987216}.

Despite decades of successful use, quantitative analysis of the convergence of Metropolis algorithms has received significant attention only in recent years \citep{pmlr-v75-dwivedi18a,pmlr-v134-chewi21a,24-AAP2058,goyal2025mixing}; for a comprehensive review, see for instance \cite{24-AAP2058}.
These convergence results quantify the mixing behavior of Metropolis algorithms, and provide theoretical guidance for tuning hyperparameters in practice.
A key quantity governing the convergence speed of Metropolis algorithms is the $L^2$-spectral gap of the associated Markov operator \citep{MR2466937,MR3155209}, and it was not until recently that \citet{24-AAP2058} derived the first explicit lower bounds on the spectral gap of RWM.
Their analysis builds on earlier works \citep{MR2178341,pmlr-v75-dwivedi18a} and develops an isoperimetric-inequality framework for estimating the spectral gap of Metropolis algorithms.
Using this framework, the authors establish explicit spectral gap bounds of RWM that achieve minimax-optimal dependence on the dimension and the condition number for smooth, strongly log-concave distributions.
A more recent work \citep{goyal2025mixing} further generalizes this framework to distributions that only satisfy a Poincar\'e or log-Sobolev inequality, thereby extends the isoperimetric approach beyond the strongly log-concave setting.

However, existing theoretical results on the convergence of Metropolis algorithms mostly focus on smooth distributions.
Such analyses typically rely on differentiability of the log-density, and often require additional regularity such as Lipschitz continuity of the gradient \citep{pmlr-v75-dwivedi18a,JMLR:v23:21-1184,24-AAP2058}.
This limits their applicability to non-smooth target distributions which are prevalent in many applications, e.g., Bayesian Lasso \citep{Park01062008,MR2370074}, image processing \citep{MR3763089,24M1694781}, and graphical model learning \citep{MAL-001}.
Such non-smoothness often arises from $\ell_1$-type regularization or ReLU-like activation functions \citep{pmlr-v97-vladimirova19a}, leading to log-densities that are not differentiable everywhere.
The non-differentiability also poses challenges for gradient-based Metropolis algorithms such as MALA, and developing effective modifications remains an active research topic \citep{MR3515019,JMLR:v23:20-527}.

In this work, we extend the isoperimetric approach to analyze the spectral gap of Metropolis algorithms for a broad class of non-smooth distributions, where the log-density can be decomposed into a smooth part and a non-smooth Lipschitz part.
We derive explicit spectral gap bounds for both RWM and MALA under appropriate step-size constraints.
For RWM, we further generalize the results of \cite{24-AAP2058} to distributions satisfying a Poincar\'e or log-Sobolev inequality by combining our analysis with \cite{goyal2025mixing}.
Our main technical contribution is to establish the close-coupling condition \citep{MR2178341} for Metropolis algorithms over non-smooth distributions.
As a by-product, we obtain a uniform lower bound on the acceptance rate of MALA under log-concavity, which plays a key role in our spectral-gap analysis. To the best of our knowledge, such a uniform bound is not available in the existing literature.
These theoretical results are further supported by numerical experiments.

The rest of the paper is organized as follows.
In \Cref{sec:Prelim}, we review Metropolis algorithms and introduce the class of non-smooth distributions considered in this work.
In \Cref{sec:Method}, we summarize the isoperimetric approach for spectral gap of Metropolis algorithms in \cite{24-AAP2058,goyal2025mixing}.
In \Cref{sec:SpecGapMH}, we present our main results on the spectral gap of RWM and MALA for non-smooth distributions.
In \Cref{sec:Numerics}, we report numerical experiments that verify the theoretical results.
We conclude in \Cref{sec:Conclusion}, and technical proofs are deferred to the appendix.

\paragraph{Notations}
For a vector $x\in\mR^d$, denote $\norme{x}_p$ as the $\ell_p$-norm of $x$, and simply write $\norme{x} = \norme{x}_2$.
Denote $f = \mcO(g)$ if there exists a universal constant $C>0$ such that $f \le C g$.
Denote $\msP(\mR^d)$ as the set of probability distributions on $\mR^d$.
Denote $\GN(m,C)$ as the Gaussian distribution with mean $m$ and covariance $C$, and denote $\GN(x;m,C)$ as its density function evaluated at $x$.
For a probability distribution $\pi$ and a function $u(x)$, denote $\E_\pi[u] = \int u(x) \pi(x) \mdd x$ as the expectation of $u$ under $\pi$.
Denote the $\pi$-weighted $L^2$-norm of $u$ as
\(
    \norm{u}_{L^2(\pi)} := \Brac{\int u^2(x) \pi(x) \mdd x}^{1/2}.
\)
Denote the weighted space $L^2(\pi) = \{ u:\mR^d\gto \mR \mid \norm{u}_{L^2(\pi)} <\infty \}$.
Let $\sfP(x,y)$ denote a Markov kernel, denote 
\(
    \sfP u(x) = \int \sfP(x,y) u(y) \mdd y,\, 
    \Brac{ \pi \otimes \sfP } (x,y) = \pi(x) \sfP(x,y).
\)

\section{Metropolis algorithms and non-smooth distributions}
\label{sec:Prelim}

\subsection{Metropolis algorithms}
\label{sec:Metropolis}

Let $\pi \in \msP(\mR^d)$ be a probability distribution on $\mR^d$ with density $\pi(x)$.
The Metropolis algorithm is a class of Markov chain Monte Carlo methods that approximately sample from $\pi$ by constructing a Markov transition kernel $\sfP$ that leaves $\pi$ invariant.
Given a proposal kernel $\sfQ$, the algorithm first proposes $x' \sim \sfQ(x,\cdot)$ and then accepts it with probability
\begin{equation}
\label{eqn:acrate}
    \alpha(x,x') = \min \left\{1,\frac{\pi(x')\sfQ(x',x)}{\pi(x)\sfQ(x,x')} \right\}.
\end{equation}
The resulting transition kernel is
\begin{equation}
\label{eqn:MHkernel}
    \sfP(x,y) = \alpha(x,y) \sfQ(x,y)
    + \Brac{ 1-\int \alpha(x,x')\sfQ(x,x')\,\mdd x'} \delta_x(y).
\end{equation}
Here $\delta_x$ is the Dirac delta measure at $x$.
Note $\sfP$ is reversible with respect to $\pi$, i.e., it satisfies the detailed balance condition 
\[  
    \pi(x) \sfP(x,y) = \pi(y) \sfP(y,x),
\]
and thus leaves $\pi$ invariant.
In this work, we focus on two Metropolis algorithms, i.e., the Random--Walk Metropolis and the Metropolis--adjusted Langevin algorithm.

\subsubsection{Random--Walk Metropolis}
\label{sec:RWM}
In RWM, the proposal is generated by a random walk step:
\begin{equation}
    x' = x + h \xi,\quad \xi\sim \GN(0,\sfI).
\end{equation}
Therefore, the proposal kernel $\sfQ_{\rm RWM}$ is given by
\begin{equation}
\label{eqn:Q_RWM}
    \sfQ_{\rm RWM}(x,x') = \GN(x';x,h^2 \sfI). 
\end{equation}
Notice since $\sfQ_{\rm RWM}$ is symmetric, the acceptance probability \eqref{eqn:acrate} reduces to
\begin{equation}
\label{eqn:acrate2}
    \alpha (x,x') = \min \left\{1,\frac{\pi(x')}{\pi(x)} \right\}.
\end{equation}

\subsubsection{Metropolis--adjusted Langevin algorithm}
\label{sec:MALA}
In MALA, the proposal is generated by a single step of the Euler-Maruyama discretization of the overdamped Langevin dynamics:
\begin{equation}
    x' = x + h \nabla \log \pi(x) + \sqrt{2h} \xi, \quad \xi \sim \GN(0,\sfI).
\end{equation}
Equivalently, the proposal kernel $\sfQ_{\rm MALA}$ is given by
\begin{equation}
\label{eqn:Q_MALA}
    \sfQ_{\rm MALA}(x,x') = \GN(x'; x + h \nabla \log \pi(x), 2h \sfI).
\end{equation}
Note that MALA is a gradient-based method, i.e., it requires the evaluation of $\nabla \log \pi(x)$ at each step.

\subsection{Non-smooth distributions}
\label{sec:NonSmthDist}

In many applications, the target distribution $\pi$ may be non-smooth.
For example, in Bayesian statistics, non-smooth priors such as the Laplace prior or spike-and-slab priors are commonly employed to induce sparsity in high-dimensional models \citep{Mitchell01121988,Park01062008}.
In machine learning, non-smooth activation functions such as ReLU can likewise lead to non-smooth posterior distributions \citep{pmlr-v97-vladimirova19a}.
Such non-smooth distributions generally do not admit gradients that are defined everywhere, which poses challenges for gradient-based sampling methods and complicates the convergence analysis of MCMC algorithms.
To address this challenge, \cite{MR3515019} introduced a proximal MALA algorithm for sampling from non-smooth distributions and established its geometric ergodicity.
Subsequent works \citep{JMLR:v23:20-527,liang2022proximal,liang2023a} proposed and analyzed various Metropolis-adjusted proximal variants.
In \cite{24M1694781}, the authors proposed a MALA-within-Gibbs algorithm, based on a different smoothing strategy, for sampling from non-smooth posteriors arising in image deblurring problems.

In this work, we consider a class of non-smooth distributions, where the log density can be decomposed into a smooth part and a non-smooth globally Lipschitz part, which is referred to as a non-smooth composite potential in \cite{JMLR:v23:20-527}.
Specifically, we make the following assumption on the target distribution $\pi$.

\begin{asm}
\label{asm:smth}
The log-density of the distribution $\pi \in\msP(\mR^d)$ can be decomposed as
\begin{equation}
\label{eqn:decomp}
    \log \pi(x) = f(x) + g(x),
\end{equation}
and the following conditions hold:
\begin{enumerate}[1.]
    \item $f(x)$ is $\sfM$-smooth, i.e.,
    \begin{equation}
    \label{eqn:Msmooth}
        \forall x,y\in\mR^d, \quad \norme{ \nabla f(x) - \nabla f(y) } \le \sfM \norme{x-y}.
    \end{equation}
    \item $ g(x)$ is $\sfL$-Lipschitz, i.e.,
    \begin{equation}
    \label{eqn:LLip}
        \forall x,y\in\mR^d, \quad \norme{ g(x) - g(y) } \le \sfL \norme{x-y}.
    \end{equation}
\end{enumerate}
\end{asm}

The decomposition \eqref{eqn:decomp} arises naturally in many applications.
In Bayesian statistics, $\pi$ often represents a posterior distribution, where $f(x)$ typically corresponds to the log-likelihood and $g(x)$ to the log-prior density.
In many sparsity-promoting problems, the prior is chosen to be of $\ell_1$-type, such as the Laplace prior, which use a non-smooth log-density.
Note the assumption is not limited to this case, and also allows situations in which the non-smooth component appears in the likelihood \citep{MR1861390}.

\Cref{asm:smth} essentially requires that the non-smooth component $g$ is globally Lipschitz, which is satisfied if the non-smoothness comes from $\ell_1$-regularization or ReLU-like potentials. 
For instance, consider the Bayesian Lasso model
\begin{equation}
    \pi(x) \propto \exp \Brac{ - \frac{1}{2}\norme{Ax-b}^2 - \mu \norme{x}_1 }. 
\end{equation}
The non-smooth part $g(x) = -\mu \norme{x}_1$ is globally $\mu$-Lipschitz.
This assumption relaxes the smoothness conditions imposed in the existing literature on spectral gap of Metropolis algorithms \citep{pmlr-v75-dwivedi18a,JMLR:v23:21-1184,24-AAP2058}, i.e., there exists $\sfM>0$ such that $ \norme{ \nabla \log \pi (x) - \nabla \log \pi(y) } \le \sfM \norme{x-y}$, which is not satisfied for non-smooth distributions.

\begin{rem}
An important class of non-smooth targets not covered in this work is the distributions with constrained support, where the density includes an indicator term and is therefore not continuous at the boundary.
There is a vast literature on MCMC for such constrained distributions; see for instance \cite{MR1733749}.
Another important class not covered here is the distributions with point masses, such as spike-and-slab priors \citep{MR2370074}, which requires a different analysis and is outside the scope of this paper.
\end{rem}

\subsection{Metropolis algorithms for non-smooth distributions}
\label{sec:MH4NonSmth}

Metropolis algorithms can still be applied to sample from non-smooth distributions, with possible minor modifications.
For RWM, no modification is needed.
For MALA, we need to specify the gradient $\nabla \log \pi(x)$ in the proposal distribution \eqref{eqn:Q_MALA}.
We follow the subgradient approach \citep{JMLR:v20:18-173}, where the gradient is replaced by the subgradient $\partial \log \pi(x)$.

Suppose $\pi$ is log-concave, i.e., $-\log \pi(x)$ is convex.
Define its subgradient as
\begin{equation}
\label{eqn:def_subdiff}
    \partial \log \pi(x) = \left\{ v \in \mR^d : \forall y\in\mR^d,~ \log \pi(y) \le \log \pi(x) + v \cdot (y-x) \right\}.
\end{equation}
It is well known that $\partial \log \pi(x)$ is non-empty for all $x\in\mR^d$ when $\pi$ is log-concave, and at the points where $\log \pi(x)$ is differentiable, $\partial \log \pi(x)$ reduces to the singleton $\{ \nabla \log \pi(x) \}$ \citep{MR2142598}.
Since $\partial \log \pi(x)$ is a set, we need to choose a specific subgradient in the proposal distribution \eqref{eqn:Q_MALA}.
A possible choice is the subgradient with minimal norm, i.e.,
\begin{equation}
    v(x) := \argmin \left\{ \norme{v} : v \in \partial \log \pi(x) \right\}.
\end{equation}
It can be easily verified that under \Cref{asm:smth}, for any choice of $v(x) \in \partial \log \pi(x)$, we can decompose
\begin{equation}
\label{eqn:v_decomp}
    v(x) = \nabla f (x) + v_{\rm s}(x), \quad v_{\rm s}(x) \in \partial g(x). 
\end{equation}
Therefore, we can draw the MALA proposal for non-smooth distributions as
\begin{equation}
    x' = x + h v(x) + \sqrt{2h} \xi, \quad \xi \sim \GN(0,\sfI).
\end{equation}
That is, $\sfQ_{\rm MALA}(x,x') = \GN(x'; x + h v(x), 2h \sfI)$. This is also known as the MSALA algorithm in \cite{ning2025metropolis}.

\section{Isoperimetric approach to the spectral gap}
\label{sec:Method}

\subsection{Spectral gap}
\label{sec:SpecGap}

The spectral gap of a Markov kernel $\sfP$ plays a crucial role in quantifying the convergence speed of the associated Markov chain to its invariant distribution $\pi$.

\begin{defn}
\label{defn:SpecGap}
Let $\sfP$ be a Markov kernel with invariant distribution $\pi$.
The spectral gap of $\sfP$ is defined as
\begin{equation}
\label{eqn:gap}
    \Gap(\sfP) = \sup \left\{ 1 - \frac{ \E_\pi \Rectbrac{ u \cdot \sfP u } }{ \E_\pi \Rectbrac{ u^2 } } : 0 \neq u \in L^2(\pi),~ \E_\pi [u] = 0 \right\}.
\end{equation}
\end{defn}

When $\sfP$ is reversible with respect to $\pi$, the spectral gap quantifies the $L^2(\pi)$-convergence rate of the associated Markov chain:
\begin{equation}
\label{eqn:L2conv}
    \forall u \in L^2(\pi), \quad \norm{\sfP^n u - \E_\pi[u] }_{L^2(\pi)} \le (1-\Gap(\sfP))^n \norm{u - \E_\pi[u] }_{L^2(\pi)}.
\end{equation}
We refer to Chapter 12 in \cite{MR2466937} for more details on spectral gap and its connections to convergence analysis of Markov chains.

A fundamental tool for estimating the spectral gap is \emph{Cheeger's inequality} \citep{MR402831}, which relates the spectral gap to the \emph{conductance} of $\sfP$:
\begin{equation}
\label{eqn:Cond}
    \Phi(\sfP) = \inf \left\{ \, \frac{ \Brac{\pi \otimes \sfP } (S \times S\matc)}{\pi(S)} \, :\; S \subseteq \mR^d,~ 0 < \pi(S) \le \frac{1}{2} \right\}.
\end{equation}

\begin{thm}[\cite{MR930082}]
\label{thm:ChIneq}
Let $\sfP$ be a reversible Markov kernel. Then it holds
\begin{equation}
\label{eqn:ChIneq}
    \frac{\kappa}{2} \Phi(\sfP)^2 \le \Gap(\sfP) \le \Phi(\sfP),
\end{equation}
where $\kappa \ge 1$ is a universal constant.
\end{thm}

By \Cref{thm:ChIneq}, it suffices to lower bound the Cheeger constant $\Phi(\sfP)$ to obtain a lower bound on the spectral gap.
There are various techniques to lower bound the Cheeger constant, among which the isoperimetric approach is particularly effective for strongly log-concave distributions \citep{pmlr-v75-dwivedi18a,24-AAP2058}.
More recently, this approach has been extended beyond log-concavity to distributions satisfying Poincar\'e or log-Sobolev inequalities \citep{goyal2025mixing}.
In the next section, we briefly review the isoperimetric approach and present a general spectral gap result derived from this framework.

\subsection{Isoperimetric approach}
\label{sec:IsoApproach}

The isoperimetric approach \citep{MR2178341,pmlr-v75-dwivedi18a,24-AAP2058} establishes a lower bound on the Cheeger constant through two key ingredients:

\begin{itemize}
    \item \emph{Three-set isoperimetric inequality} for $\pi$.
    There exists an increasing function \( F:(0,\frac{1}{2} ] \to \mR_+ \) such that, for any measurable decomposition
    \( \mR^d = S_1 \sqcup S_2 \sqcup S_3\),
    it holds that
    \begin{equation}
    \label{eqn:3setIsoIneq}
        \pi(S_3) \ge d(S_1,S_2) \, F \Brac{ \min \{ \pi(S_1), \pi(S_2) \} },
    \end{equation}
    where $d(S_1,S_2) := \inf_{x\in S_1,\,y\in S_2} \norme{ x-y }$.
    \item \emph{Close-coupling condition} for \(\sfP\).
    There exist constants $\epsilon,\delta>0$ such that
    \begin{equation}
    \label{eqn:CloseCouple}
        \norme{ x-y } < \delta ~\St~ \TV (\delta_x \sfP,\delta_y \sfP) < 1-\epsilon .
    \end{equation}
    Here $\TV$ denotes the total variation distance between two probability measures.
    \begin{equation}
    \label{eqn:TV}
        \TV(\mu,\nu) = \sup_{A\subseteq\mR^d} |\mu(A) - \nu(A)| = \frac{1}{2} \int | \mu(x) - \nu(x) | \mdd x.
    \end{equation}
\end{itemize}

We make several remarks on the three-set isoperimetric inequality \eqref{eqn:3setIsoIneq}.
First, it is worth noting that \eqref{eqn:3setIsoIneq} depends only on the target distribution $\pi$ and is independent of the Markov kernel \(\sfP\).
Second, \eqref{eqn:3setIsoIneq} is equivalent to that $F$ is a minorant of the \emph{isoperimetric profile} $I_\pi$ of $\pi$; see Section 3 in \cite{24-AAP2058} for details.
\cite{24-AAP2058} also show several lower bounds of $I_\pi$ for log-concave distributions, yielding explicit choices of $F$.
Finally, \cite{goyal2025mixing} recently generalize \eqref{eqn:3setIsoIneq} to the form
\begin{equation}
\label{eqn:3setIsoIneq_gen}
    \pi(S_3) \ge \Upsilon \Brac{ d(S_1,S_2) }
    F \Brac{\min\{\pi(S_1),\pi(S_2)\}},
\end{equation}
where $\Upsilon:\mR_+\to\mR_+$ is an increasing function.
Note \eqref{eqn:3setIsoIneq_gen} reduces to \eqref{eqn:3setIsoIneq} when $\Upsilon = \id$.
This extension applies to a broader class of target distributions, including those satisfying Poincar\'e or log-Sobolev inequalities.

Now we present a general spectral gap result based on the isoperimetric approach.

\begin{prop}
\label{prop:SpecGap_via_Iso}
Let $\sfP$ be a reversible Markov kernel with invariant distribution $\pi \in \msP(\mR^d)$.
Suppose the generalized three-set isoperimetric inequality \eqref{eqn:3setIsoIneq_gen} holds for $\pi$ with some increasing functions $F:(0,\frac{1}{2}] \to \mR_+$ and $\Upsilon:\mR_+ \to \mR_+$, and the close-coupling condition \eqref{eqn:CloseCouple} holds for $\sfP$ with some $\epsilon,\delta>0$.
Then the spectral gap of $\sfP$ admits the lower bound
\begin{equation}
\label{eqn:SpecGap_via_Iso}
    \Gap(\sfP) \ge \frac{\epsilon^2}{8} \Rectbrac{ \sup_{\theta\in[0,1]} \min \Big\{ 1-\theta , \Upsilon(\delta) \inf_{t\in(0,\frac{1}{2}]} \frac{ F(\theta t) }{ 2t } \Big\} }^2.
\end{equation}
When $F$ is concave on $(0,\frac{1}{2}]$, we can obtain a simplified lower bound
\begin{equation}
\label{eqn:Gap_bound}
    \Gap(\sfP) \ge \frac{\epsilon^2}{8} \Rectbrac{ \min \Big\{ \frac{1}{2} , \Upsilon(\delta) F(\tfrac{1}{4}) \Big\} }^2.
\end{equation}
\end{prop}

\Cref{prop:SpecGap_via_Iso} is a slight modification of Proposition 3.11 in \cite{goyal2025mixing}, and we provide a self-contained proof in \Cref{app:pf_SpecGap_via_Iso} for completeness.

\subsection{Three-set isoperimetric inequality}
\label{sec:3setIsoIneq}

The three-set isoperimetric inequality relies on certain convexity-type properties of the target distribution \citep{MR2507637,24-AAP2058,goyal2025mixing}.
We summarize below several results from \cite{24-AAP2058} and \cite{goyal2025mixing}, which establish \eqref{eqn:3setIsoIneq} or \eqref{eqn:3setIsoIneq_gen} under assumptions such as strongly log-concavity, and Poincar\'e or log-Sobolev inequalities for the target distribution.

\begin{defn}
\label{defn:convexity}
Let $\pi \in \msP(\mR^d)$ be a probability distribution.

\begin{itemize}
    \item $\pi$ is called \textbf{$\sfm$-strongly log-concave} if its log-density $\log \pi$ is $\sfm$-strongly concave, i.e., for all $x,y\in\mR^d$ and $\lambda\in(0,1)$,
    \begin{equation}
    \label{eqn:convex}
        \log \pi ( \lambda x + (1-\lambda) y ) \ge \lambda \log \pi(x) + (1-\lambda) \log \pi(y) + \frac{\sfm}{2} \lambda (1-\lambda) \norme{x-y}^2.
    \end{equation}

    \item $\pi$ is said to satisfy a \textbf{log-Sobolev inequality} with constant $\CLSI>0$ if for all locally Lipschitz functions $u:\mR^d\to\mR$,
    \begin{equation}
    \label{eqn:logSobolev}
        \E_\pi \Rectbrac{ |u|^2 \log \frac{ |u|^2 }{ \E_\pi [|u|^2] } } \le \CLSI \; \E_\pi \Rectbrac{ \norme{ \nabla u }^2 }.
    \end{equation}

    \item $\pi$ is said to satisfy a \textbf{Poincar\'e inequality} with constant $\CPI>0$ if for all locally Lipschitz functions $u:\mR^d\to\mR$,
    \begin{equation}
    \label{eqn:Poincare}
        \E_\pi \Rectbrac{ \norme{ u - \E_\pi [u] }^2 } \le \CPI \; \E_\pi \Rectbrac{ \norme{ \nabla u }^2 }.
    \end{equation}
\end{itemize}
\end{defn}

It is well known \citep{MR3155209} that the above convexity-type properties satisfy:
\begin{center}
    strongly log-concavity $~\St~$ log-Sobolev inequality $~\St~$ Poincar\'e inequality.
\end{center}
Under strong log-concavity, the following proposition holds.
\begin{prop}[\cite{24-AAP2058}, Lemma 27]
\label{prop:convex_IsoIneq}
Suppose $\pi \in \msP(\mR^d)$ is $\sfm$-strongly log-concave.
Then the three-set isoperimetric inequality \eqref{eqn:3setIsoIneq} holds with
\begin{equation}
\label{eqn:F_convex}
    F(t) = t \Brac{ \frac{2 \sfm \log t^{-1} }{\pi \log 2} }^{1/2}.
\end{equation}
\end{prop}

The strongly log-concavity assumption can be relaxed to log-Sobolev or Poincar\'e inequalities, yielding the following results.
\begin{prop}[\cite{goyal2025mixing}, Propositions 4.1 and 4.3]
\label{prop:functional_IsoIneq}
Let $\pi \in \msP(\mR^d)$.
Assume that $\pi$ satisfies either a log-Sobolev inequality or a Poincar\'e inequality, then the generalized three-set isoperimetric inequality \eqref{eqn:3setIsoIneq_gen} holds, and the functions $F$ and $\Upsilon$ can be chosen as follows:
\begin{itemize}
    \item \textbf{log-Sobolev case:} if $\pi$ satisfies a log-Sobolev inequality with constant $\CLSI>0$, then
    \begin{equation}
    \label{eqn:F_Ups_LSI}
        F(t) = \frac{t}{2}\log\frac{2}{t}, \quad \Upsilon(\delta) = \frac{\delta^2}{\CLSI+(1+\mee^{-1})\delta^2}.
    \end{equation}
    \item \textbf{Poincar\'e case:} if $\pi$ satisfies a Poincar\'e inequality with constant $\CPI>0$, then
    \begin{equation}
    \label{eqn:F_Ups_PI}
        F(t) = t, \quad \Upsilon(\delta) = \frac{\delta^2}{16(\CPI+4\delta^2)}.
    \end{equation}
\end{itemize}
\end{prop}

\begin{rem}
In \cite{goyal2025mixing}, the authors consider general $L^p$-log Sobolev inequalities and $L^p$-Poincar\'e inequalities for $p\ge 1$.
Here we consider only the usual case $p=2$.
Also note the log-Sobolev and Poincar\'e constants are defined as the reciprocals of those in \cite{goyal2025mixing}, which leads to slight differences in the expressions of $\Upsilon(\delta)$.
\end{rem}

\section{Spectral gap of Metropolis algorithms for non-smooth distributions}
\label{sec:SpecGapMH}

\subsection{Spectral gap of Random--Walk Metropolis}
\label{sec:SpecGapRWM}

\subsubsection{Close coupling condition}
\label{sec:CloseCoupleRWM}

By the discussions in \Cref{sec:IsoApproach}, to establish a spectral gap for Random--Walk Metropolis, it remains to verify the close coupling condition \eqref{eqn:CloseCouple}.
Previous work \citep{24-AAP2058} relies on the smoothness on the target distribution, i.e., the log-density $\log \pi$ has bounded Hessian.
To handle the non-smoothness, we develop new ways to establish the close coupling condition, aiming to obtain explicit and sharper dependence of the step size $h$ on the parameters $d$, $\sfM$, $\sfL$.

\begin{lem}
\label{lem:CloseCoupleRWM}
Let $\pi \in \msP(\mR^d)$ satisfy \Cref{asm:smth}.
Consider the Random--Walk Metropolis kernel $\sfP$ with step size $h$.
Assume that
\begin{equation}
\label{eqn:h_RWM}
    0 < h \le \frac{1}{\sqrt{d}}
    \min\left\{
        \frac{1}{2\sqrt{\sfM}},
        \frac{1}{16\sfL}
    \right\}.
\end{equation}
Then $\sfP$ satisfies a close coupling condition:
\begin{equation}
\label{eqn:CloseCoupleRWM}
    \norme{x-y} < \frac{h}{4} ~\St~ \TV(\delta_x \sfP,\delta_y \sfP) < \frac{3}{4}.
\end{equation}
\end{lem}

Proof of \Cref{lem:CloseCoupleRWM} is provided in \Cref{app:pf_CloseCoupleRWM}.

\begin{rem}
\label{rem:CloseCoupleRWM}
There are two constraints on the step size: \(h \lesssim d^{-1/2} \sfM^{-1/2} \) and \(h \lesssim d^{-1/2}\sfL^{-1}\).
The first constraint is known to be necessary to obtain a non-degenerate spectral gap for smooth target distributions; see \cite{24-AAP2058}.
The second constraint is new and arises from the presence of the non-smooth component.
We comment that the constants in \eqref{eqn:h_RWM} are not optimized, and may be improved by more careful analysis. But the scaling with respect to $d$, $\sfM$ and $\sfL$ is tight.
A qualitative interpretation of the results is that there exisits a step size regime $ 0 < h \le \mean{h}$, where $\mean{h} \asymp d^{-1/2} \min\{\sfM^{-1/2}, \sfL^{-1}\}$, such that the close coupling condition holds.

\end{rem}

\subsubsection{Spectral gap results}
\label{sec:RWM_SpecGap}

We are now ready to present our main spectral gap results for Random--Walk Metropolis under non-smooth distributions.
Under the strongly log-concave condition, we have

\begin{thm}
\label{thm:SpecRWN_convex}
Let $\pi\in\msP(\mR^d)$ satisfy \Cref{asm:smth}.
Assume further that $\pi$ is $\sfm$-strongly log-concave with $\sfm \le \sfM$.
Consider the Random--Walk Metropolis kernel $\sfP$ with step size $h$ satisfying \eqref{eqn:h_RWM}.
Then the spectral gap of $\sfP$ admits the lower bound
\begin{equation}
    \Gap (\sfP) \ge \sfC\sfm h^2.
\end{equation}
Here $\sfC$ is some universal constant.
\end{thm}

\begin{rem}
\label{rem:SpecRWM_convex}
(1) Take $h = \frac{h_0}{\sqrt{d}} \min \left\{ \frac{1}{2\sqrt{\sfM}}, \frac{1}{16\sfL} \right\} $ where $0<h_0\le 1$, we get 
\begin{equation}
    \Gap(\sfP) \ge \sfC'h_0^2 \cdot \min \left\{ d^{-1} \kappa^{-1} , d^{-1} \sfm \sfL^{-2} \right\}.
\end{equation}
Here we denote $\kappa = \sfM/\sfm$ as the condition number.
The scaling $\mcO(d^{-1} \kappa^{-1})$ is known to be optimal for smooth strongly log-concave distributions \citep{24-AAP2058}.
The scaling $\mcO(d^{-1} \sfm \sfL^{-2})$ comes from the additional constraint due to the non-smoothness. 

\noindent
(2) The universal constant $\sfC$ obtained in our proof is $\frac{1}{2^{13}\pi}$, which is clearly not sharp due to several relaxations in the argument. We do not attempt to optimize it in this work.
A qualitative interpretation of the result is that there exists a step size regime $ 0 < h \le \mean{h}$, where $\mean{h} \asymp d^{-1/2} \min\{\sfM^{-1/2}, \sfL^{-1}\}$, such that the spectral gap scales as $\mcO(\sfm h^2)$.

\end{rem}

\begin{proof}[Proof for \Cref{thm:SpecRWN_convex}]
By \Cref{lem:CloseCoupleRWM}, the close coupling condition \eqref{eqn:CloseCouple} holds for $\sfP$ with $\epsilon = \frac{1}{4}, \delta = \frac{h}{4}$.
By \Cref{prop:convex_IsoIneq}, the three-set isoperimetric inequality \eqref{eqn:3setIsoIneq} holds for $\pi$ with $F$ in \eqref{eqn:F_convex}.
Therefore, by \Cref{prop:SpecGap_via_Iso}, we have
\begin{equation}
    \Gap(\sfP) \ge \frac{\epsilon^2}{8} \Rectbrac{ \min \Big\{ \frac{1}{2} , \delta F(\tfrac{1}{4}) \Big\} }^2 = \frac{1}{2^7} \min \Big\{ \frac{1}{2^2} , \frac{ \sfm h^2}{2^6 \pi} \Big\} = \sfC \sfm h^2.
\end{equation}
where we take $ \sfC = \frac{1}{2^{13}\pi}  > 3.885 \times 10^{-5} $.
Note here we use $\sfm \le \sfM$.
\end{proof}

Following similar arguments, we can also establish spectral gap results under log-Sobolev or Poincar\'e inequalities.

\begin{thm}
\label{thm:SpecRWN_func}
Let $\pi\in\msP(\mR^d)$ satisfy \Cref{asm:smth}.
Assume further that $\pi$ satisfies either a log-Sobolev inequality with constant $\CLSI>0$, or a Poincar\'e inequality with constant $\CPI>0$.
Consider the Random--Walk Metropolis kernel $\sfP$ with step size $h$ satisfying \eqref{eqn:h_RWM}.
Then the spectral gap of $\sfP$ admits the lower bound
\begin{equation}
    \Gap (\sfP) \ge \sfC' \cdot
    \begin{cases}
        \Brac{ \mathsf{C_{LSI}^2} + h^4 }^{-1} h^4 , & \text{if $\pi$ satisfies a log-Sobolev inequality}, \\[5pt]
        \Brac{ \mathsf{C_{PI}^2} + h^4 }^{-1} h^4 , & \text{if $\pi$ satisfies a Poincar\'e inequality}.
    \end{cases}
\end{equation}
Here $\sfC'$ is some universal constant.
\end{thm}

\begin{rem}
Compared to \Cref{thm:SpecRWN_convex}, the spectral gap in \Cref{thm:SpecRWN_func} has a degenerate dependence on the step size $h$, i.e., scaling as $\mcO(h^{4})$ rather than $\mcO(h^{2})$.
This is due to the weaker convexity-type assumptions.
In \Cref{prop:functional_IsoIneq}, $\Upsilon(\delta) \asymp \delta^2$ as $\delta \to 0$, in contrast to the linear dependence $\Upsilon(\delta) = \delta$ in \Cref{prop:convex_IsoIneq}.
Such degeneracy in the isoperimetric inequality is unavoidable when only $L^2$-functional inequalities are assumed; see \cite{MR2507637}.
We conjecture, however, that the $\mcO(h^{4})$ scaling of the spectral gap is merely an artifact of the isoperimetric approach.
\end{rem}

\begin{proof}[Proof for \Cref{thm:SpecRWN_func}]
By \Cref{lem:CloseCoupleRWM}, the close coupling condition \eqref{eqn:CloseCouple} holds for $\sfP$ with $\epsilon = \frac{1}{4}, \delta = \frac{h}{4}$.
When $\pi$ satisfies a log-Sobolev inequality with constant $\CLSI>0$, by \Cref{prop:functional_IsoIneq} and \Cref{prop:SpecGap_via_Iso}, we have
\begin{equation}
    \Gap(\sfP) \ge \frac{\epsilon^2}{8} \Rectbrac{ \min \Big\{ \frac{1}{2} , \frac{1}{4} \frac{\delta^2}{\CLSI+(1+\mee^{-1})\delta^2} \Big\} }^2 \ge \sfC_1 \cdot \frac{ h^4 }{ \mathsf{C_{LSI}^2} + h^4 }.
\end{equation}
where we can take $\sfC_1 = 2^{-20}$.
Likewise, when $\pi$ satisfies a Poincar\'e inequality with constant $\CPI>0$, by \Cref{prop:functional_IsoIneq} and \Cref{prop:SpecGap_via_Iso}, we have
\begin{equation}
    \Gap(\sfP) \ge \frac{\epsilon^2}{8} \Rectbrac{ \min \Big\{ \frac{1}{2} , \frac{1}{8} \log 8 \cdot \frac{\delta^2}{16(\CPI+4\delta^2)} \Big\} }^2 \ge \sfC_2 \cdot \frac{ h^4 }{ \mathsf{C_{PI}^2} + h^4 },
\end{equation}
where we can take $\sfC_2 = (\log 8)^2 \cdot 2^{-23} $.
\end{proof}

\subsection{Spectral gap of Metropolis--adjusted Langevin algorithm}
\label{sec:SpecGapMALA}

\subsubsection{Close coupling condition}
\label{sec:CloseCoupleMALA}

Similarly to the random-walk Metropolis case, it suffices to establish the close-coupling condition \eqref{eqn:CloseCouple} for MALA.
However, this condition would require a \emph{uniform} lower bound on the acceptance rate, which, to the best of our knowledge, is not available in the literature.
Existing results \citep{MR3020951,pmlr-v75-dwivedi18a,JMLR:v23:21-1184} only provide lower bounds on bounded regions, which are insufficient for our purpose.
Nevertheless, under strong log-concavity assumption, we can prove a uniform lower bound on the acceptance rate, even in the presence of non-smoothness.
We need slightly refined convexity assumptions as follows.

\begin{asm}
\label{asm:MALA}
The distribution $\pi \in\msP(\mR^d)$ satisfies \Cref{asm:smth} and is $\sfm$-strongly log-concave with $\sfm \le \sfM$.
Moreover, the function $f(x)$ in the decomposition \eqref{eqn:v_decomp} is concave.
\end{asm}

\begin{lem}
\label{lem:ac_MALA}
Let $\pi\in\msP(\mR^d)$ satisfy \Cref{asm:MALA}.
Consider the MALA kernel $\sfP$ with step size $h$ such that
\begin{equation}
\label{eqn:h_MALA}
    0 < h \le \frac{1}{200} \cdot \min \left\{ \frac{1}{d\kappa \sfM} , \frac{1}{d\sfL^2} \right\}. 
\end{equation}
Here $\kappa = \sfM/\sfm$ is the condition number.
Then the acceptance rate \eqref{eqn:acrate} of MALA proposal admits a uniform lower bound
\begin{equation}
    \inf_x \E_{y \sim \sfQ(x,y)} [\alpha(x,y)] > \frac{13}{20}.
\end{equation}
\end{lem}

Proof of \Cref{lem:ac_MALA} is provided in \Cref{app:pf_ac_MALA}.
Note here $\frac{13}{20}$ is chosen for convenience and can be replaced by any constant larger than $\frac{1}{2}$, which is sufficient for our purpose.
As in the random-walk Metropolis case, the constants in \eqref{eqn:h_MALA} are not optimized, and the results can be interpreted qualitatively (see \Cref{rem:CloseCoupleRWM}).

\begin{rem}
The step size condition \eqref{eqn:h_MALA} is more restrictive than the known results for log-concave distributions, which only requires \( h \lesssim d^{-1/2} \sfM^{-1} \text{polylog}(d,\kappa) \) \citep{JMLR:v23:21-1184}.
This suboptimality comes from technical limitations in our analysis, and leads to a non-optimal spectral gap bound.
However, relaxing the step-size constraint appears non-trivial and is left for future work.
\end{rem}

Using \Cref{lem:ac_MALA}, we can establish the close coupling condition for MALA as follows.
\begin{lem}
\label{lem:CloseCoupleMALA}
Let $\pi \in \msP(\mR^d)$ satisfy \Cref{asm:MALA}.
Consider the MALA kernel $\sfP$ with step size $h$ satisfying \eqref{eqn:h_MALA}.
Then $\sfP$ satisfies a close coupling condition:
\begin{equation}
\label{eqn:CloseCoupleMALA}
    \norme{x-y} < \frac{\sqrt{2h}}{10} ~\St~ \TV(\delta_x \sfP,\delta_y \sfP) < \frac{4}{5}.
\end{equation}
\end{lem}

Proof of \Cref{lem:CloseCoupleMALA} is provided in \Cref{app:pf_CloseCoupleMALA}.

\subsubsection{Spectral gap result}
\label{sec:MALA_SpecGap}
Now we present our main spectral gap result for MALA under non-smooth strongly log-concave distributions.
\begin{thm}
\label{thm:SpecMALA_convex}
Let $\pi\in\msP(\mR^d)$ satisfy \Cref{asm:MALA}.
Consider the MALA kernel $\sfP$ with step size $h$ satisfying \eqref{eqn:h_MALA}.
Then the spectral gap of $\sfP$ admits the lower bound
\begin{equation}
    \Gap (\sfP) \ge \sfC'' \sfm h.
\end{equation}
Here $\sfC''$ is some universal constant.
\end{thm}

\begin{rem}
Take $h = \frac{h_0}{200} \cdot \min \left\{ \frac{1}{d \kappa \sfM} , \frac{1}{d \sfL^2} \right\} $ where $0<h_0\leq 1$, we get
\[
    \Gap(\sfP) \geq \sfC'h_0 \cdot \min \left\{ d^{-1} \kappa^{-2}, d^{-1} \sfm \sfL^{-2} \right\}.
\]
The scaling \(\mcO(d^{-1}\kappa^{-2})\) is suboptimal and arises from technical artifacts in the step-size constraint.
In contrast, for smooth strongly log-concave targets, \cite{JMLR:v23:21-1184} shows that the minimax mixing time of MALA is
\(\mcO\bigl(\kappa\sqrt{d}\,\mathrm{polylog}(d,\kappa)\bigr)\),
or equivalently, that the spectral gap scales as
\(\mcO\bigl(d^{-1/2}\kappa^{-1}\mathrm{polylog}(d,\kappa)\bigr)\).
The term \(d^{-1}\sfm\,\sfL^{-2}\) arises from the additional constraint induced by the non-smooth component, as in the random-walk Metropolis case.
The constant $\sfC''$ obtained in our proof is $\frac{1}{40000 \pi}$, which is far from sharp, and the obtained lower bound should be interpreted qualitatively (see \Cref{rem:SpecRWM_convex}).
\end{rem}

\begin{proof}[Proof for \Cref{thm:SpecMALA_convex}]
By \Cref{lem:CloseCoupleMALA}, the close coupling condition \eqref{eqn:CloseCouple} holds with 
\(
    \epsilon = \frac{1}{5}, \delta = \frac{\sqrt{2h}}{10}.
\)
By \Cref{prop:convex_IsoIneq}, the three-set isoperimetric inequality \eqref{eqn:3setIsoIneq} holds for $\pi$ with $F$ in \eqref{eqn:F_convex}.
Therefore, by \Cref{prop:SpecGap_via_Iso}, we have
\begin{equation}
    \Gap(\sfP) \ge \frac{\epsilon^2}{8} \Rectbrac{ \min \Big\{ \frac{1}{2} , \delta F(\tfrac{1}{4}) \Big\} }^2 = \frac{1}{200} \min \Big\{ \frac{1}{4} , \frac{ \sfm h}{200 \pi} \Big\} = \sfC'' \sfm h.
\end{equation}
where we can take \( \sfC'' = \frac{1}{40000 \pi} > 7.957 \times 10^{-6} \).
\end{proof}

\section{Numerical experiments}
\label{sec:Numerics}

\subsection{Bayesian Lasso}
\label{sec:Num_BayeLasso}
We first consider a simple Bayesian Lasso model \citep{Park01062008}.
Let $A \in \mR^{d \times d}$ be a design matrix, and $y \in \mR^d$ be an observation vector.
The posterior distribution of the regression coefficient $x \in \mR^d$ is given by
\begin{equation}
\label{eqn:BayeLasso}
    \pi(x) \propto \exp \Brac{ - \frac{1}{2} \norme{Ax - y}^2 - \lambda \norme{x}_1 }.
\end{equation}
Here $\lambda>0$ is a regularization parameter, which controls the sparsity of the solution.
The design matrix $A$ is generated via $A = \Sigma V\matT$, where $\Sigma = \diag(\sigma_1, \dots, \sigma_d)$, $\sigma_i = i^{-\alpha}$ for some $\alpha > 0$, and $V \in \mR^{d \times d}$ is a random orthogonal matrix.
By construction, the condition number of $\pi$ is $\kappa = d^{2\alpha}$.
In the experiments, we fix $d=10$ and $\alpha=0.5$.
It is easy to verify that $\pi$ satisfies \Cref{asm:smth} with $f(x) = - \frac{1}{2} \norme{Ax - y}^2$, $g(x) = - \lambda \norme{x}_1$, $\sfM = \sigma_1^2$ and $\sfL = \lambda$.

We implement the RWM and MALA algorithms to sample from $\pi$, and estimate the spectral gap of the associated Metropolis kernel from the sampled chain via the autocorrelation method \citep{10.1214/ss/1177011137}.
Given a stationary chain $\{X_k\}_{k=0}^{K}$, the integrated autocorrelation time (IACT) of a test function $f:\mR^d \to \mR$ is defined as
\begin{equation}
    \mathrm{IACT} (f) = 1 + 2 \sum_{k=1}^\infty \rho_f(k), \quad \rho_f(k) := \frac{ \Cov( f(X_0), f(X_k) ) }{ \Var_\pi(f) }.
\end{equation}
Note that the spectral gap can be related to the IACT via \citep{10.1214/ss/1177011137}
\begin{equation}
    \Gap = \frac{2}{1+\sup_f \mathrm{IACT}(f)}.
\end{equation}
In practice, we approximate $\sup_f \mathrm{IACT}(f)$ by taking the maximum over a selected class of test functions
\(
    f \in \mcF := \{ \log \pi(x), \norme{x}_1, \norme{x}^2, \xi\matT x, \xi\matT \cos (x) \},
\)
where $\xi \in \mR^d$ is a standard Gaussian vector and $\cos (x)$ is applied elementwise.
Here the test functions are chosen empirically to ensure good estimation of the spectral gap, and we find that they work well in our experiments.
We then estimate the spectral gap by
\begin{equation}
    \widehat{\Gap} = 2 \Rectbrac{1+\max_{f\in\mathcal{F}} \; \widehat{\mathrm{IACT}}(f)}^{-1},
\end{equation}
where $\widehat{\mathrm{IACT}}(f)$ is the empirical IACT computed from the sampled chain.
This approach provides accurate estimates in our experiments.

\begin{figure}[htbp]
    \centering
    \vspace{10pt}
    \includegraphics[width=0.45\textwidth]{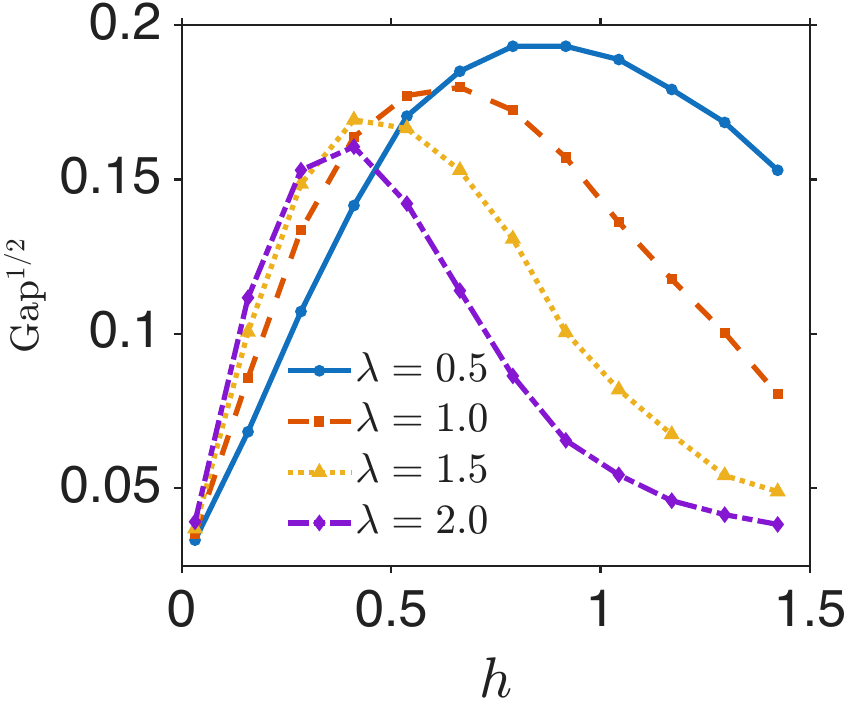} \qquad
    \includegraphics[width=0.45\textwidth]{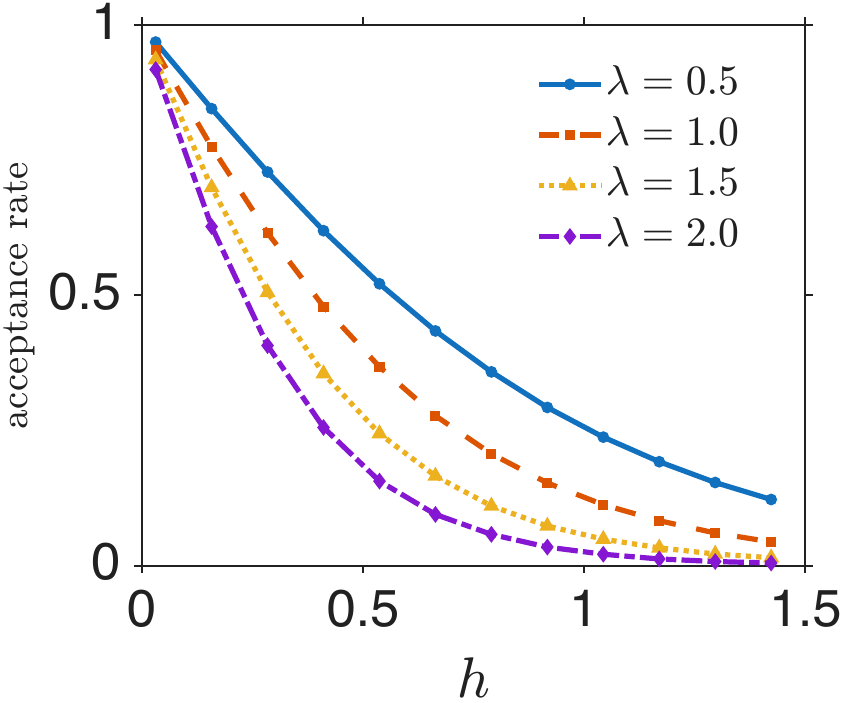}
    \vspace{-5pt}
    \caption{Estimated square root of the spectral gap (left) and average acceptance rate (right) of RWM for Bayesian Lasso model under different step size $h$ and regularization parameter $\lambda$}
    \label{fig:RWM_BayeLasso}
\end{figure}

\begin{figure}[htbp]
    \centering
    \vspace{10pt}
    \includegraphics[width=0.45\textwidth]{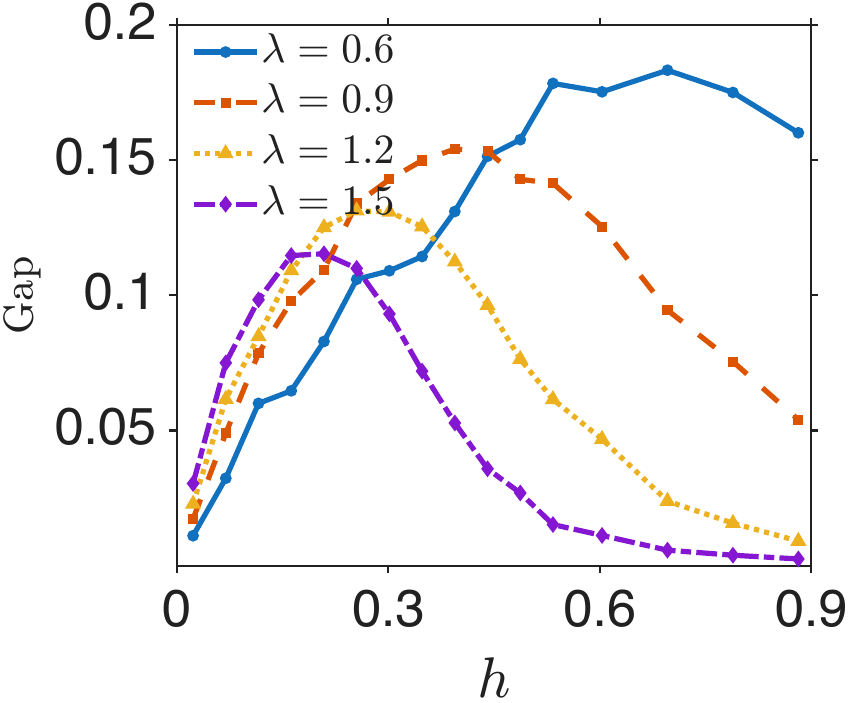} \qquad
    \includegraphics[width=0.45\textwidth]{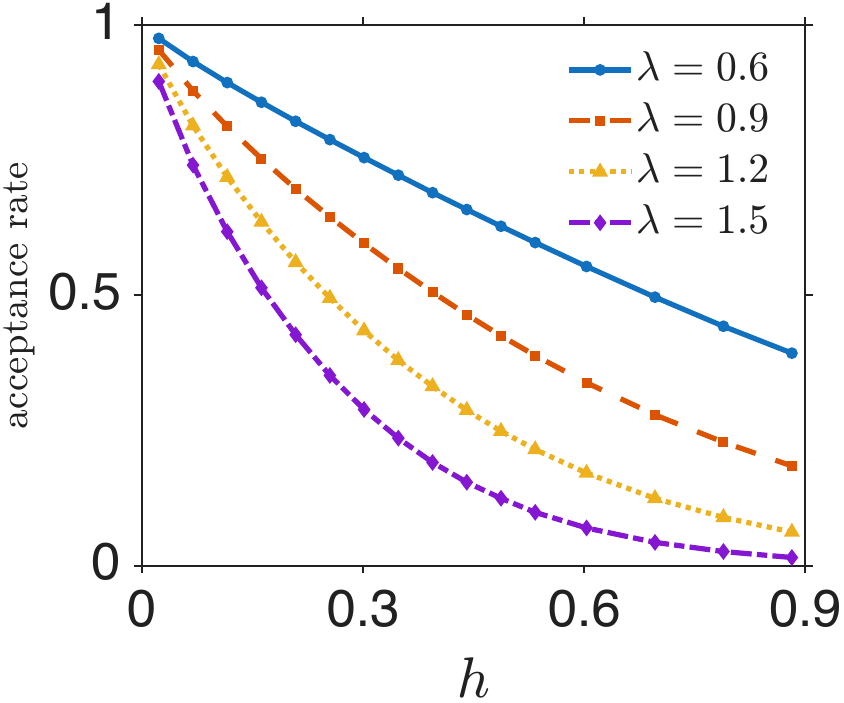}
    \vspace{-5pt}
    \caption{Estimated spectral gap (left) and average acceptance rate (right) of MALA for Bayesian Lasso model under different step size $h$ and regularization parameter $\lambda$}
    \label{fig:MALA_BayeLasso}
\end{figure}

In the experiments, we test under different step size $h$ and regularization parameter $\lambda$.
For each setting, we run $40$ independent trials with chain length $K = 10^5$, preceded by an additional $10^4$ burn-in steps, and report averages over trials.
The estimated spectral gaps and average acceptance rates are shown in \Cref{fig:RWM_BayeLasso,fig:MALA_BayeLasso} for RWM and MALA, respectively.

From \Cref{fig:RWM_BayeLasso}, we observe that the estimated spectral gap of RWM scales as $\mcO(h^2)$ for small $h$, which agrees with our theoretical result in \Cref{thm:SpecRWN_convex}.
Note the $y$-axis in the left figure is in the square root scale to better illustrate the $\mcO(h^2)$ scaling.
Moreover, the $\mcO(h^2)$ scaling region becomes smaller as $\lambda$ increases (which corresponds to stronger non-smoothness), qualitatively verifying the step size constraint due to the non-smooth component in \eqref{eqn:h_RWM}.
This can be partially explained by the fact that larger $\lambda$ leads to smaller acceptance rates, as shown in the right figure.

In \Cref{fig:MALA_BayeLasso}, we observe that the estimated spectral gap of MALA scales as $\mcO(h)$ for small $h$, which agrees with our theoretical result in \Cref{thm:SpecMALA_convex}.
Similar to RWM, the $\mcO(h)$ scaling region becomes smaller as $\lambda$ increases, qualitatively verifying the step size constraint in \eqref{eqn:h_MALA}.
The estimation of the spectral gap is less stable compared to RWM, possibly due to the empirical approximation error of IACT and suboptimality of the test functions.

Finally, we emphasize that the numerical experiments are designed only to \emph{qualitatively} validate the theoretical results, e.g., the predicted scaling of the spectral gap with respect to the step size $h$, and the step size constraints with respect to the regularization parameter $\lambda$.
Because the constants in the theoretical bounds are not sharp, exact quantitative agreement is not expected. Nonetheless, the behaviors observed are consistent with our theory.

\subsection{Logistic regression with Laplace prior}
\label{sec:Num_LogReg}
Consider a Bayesian logistic regression model with a non-smooth Laplace prior \citep{MR2227368,MR2370074}.
Given a dataset $\{(a_i,y_i)\}_{i=1}^N$ where $a_i \in \mR^d$ is a feature vector, and $y_i \in \{0,1\}$ is a binary label, the posterior distribution of the regression coefficient $x \in \mR^d$ is given by
\begin{equation}
\label{eqn:LogReg}
    \pi(x) \propto \exp \Brac{ - \sum_{i=1}^N \Brac{ \log ( 1 + \exp (a_i\matT x) ) - y_i a_i\matT x } - \lambda \norme{x}_1 }.
\end{equation}
Here $\lambda>0$ is the regularization parameter.
Note that the negative log-likelihood is convex, but not globally strongly convex. This model therefore provides a representative test case for our theory under Poincar\'e inequalities.
(since log-concave distributions satisfy Poincar\'e inequalities \citep{MR2507637})
Nevertheless, we note that it behaves like a strong log-concave distribution in a local region around its minimum.

We generate synthetic data as follows.
The feature vectors $\{a_i\}_{i=1}^N$ are sampled i.i.d. from $\GN(0,\sfI)$.
The true coefficient $x_* \in \mR^d$ is sparse, with $10\%$ non-zero entries, whose locations are chosen uniformly at random and whose values are drawn from $\GN(0,1)$.
The binary labels $\{y_i\}_{i=1}^N$ are generated via $y_i \sim \text{Bernoulli}(\sigma(a_i\matT x_*))$, where $\sigma(t) = \frac{1}{1+\mee^{-t}}$ is the sigmoid function.

\begin{figure}[htbp]
    \centering
    \vspace{10pt}
    \includegraphics[width=0.45\textwidth]{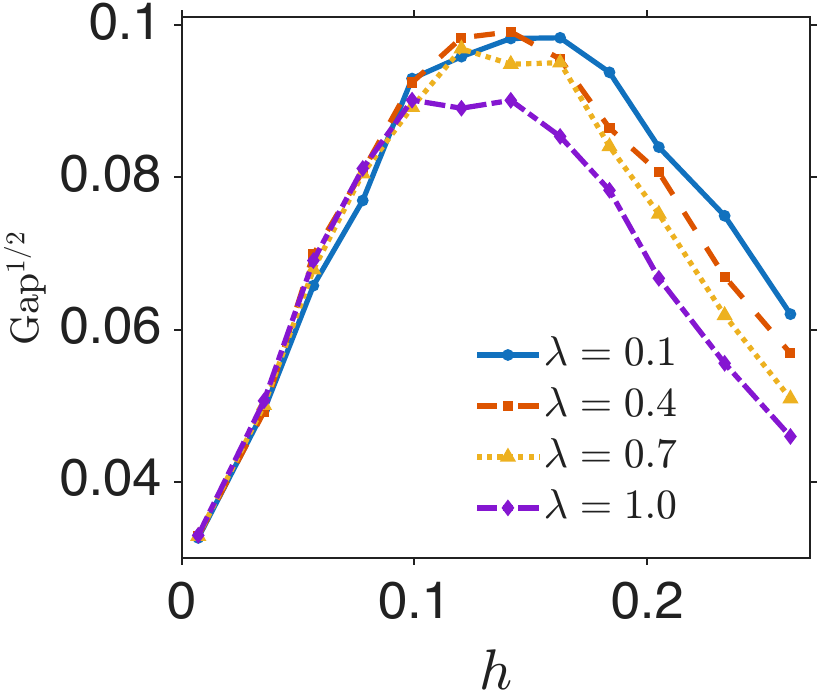} \qquad
    \includegraphics[width=0.45\textwidth]{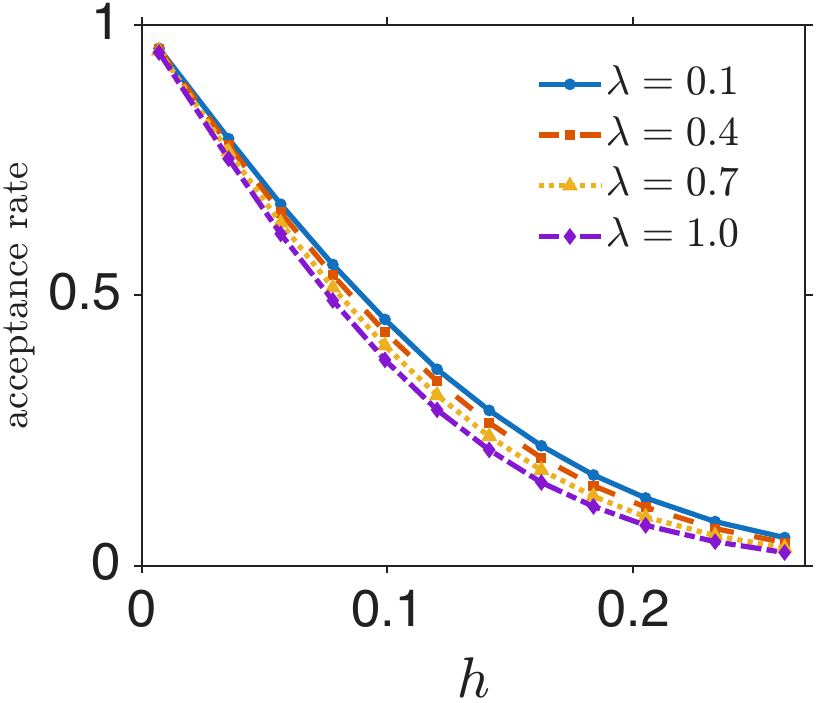}
    \vspace{-5pt}
    \caption{Estimated square root of the spectral gap (left) and average acceptance rate (right) of RWM for Bayesian Logistic regression model under different step size $h$ and regularization parameter $\lambda$}
    \label{fig:RWM_LogReg}
\end{figure}

\begin{figure}[htbp]
    \centering
    \vspace{10pt}
    \includegraphics[width=0.45\textwidth]{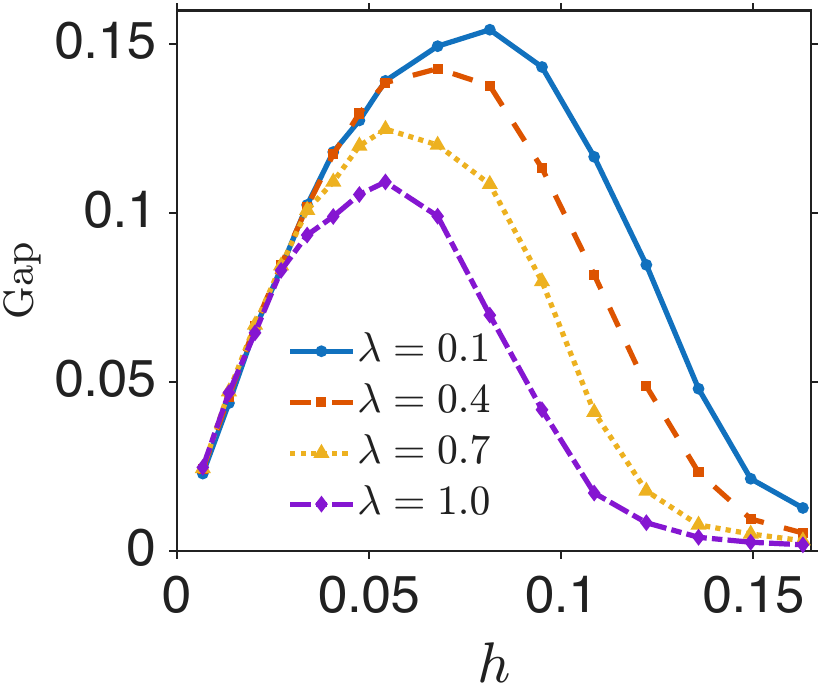} \qquad
    \includegraphics[width=0.45\textwidth]{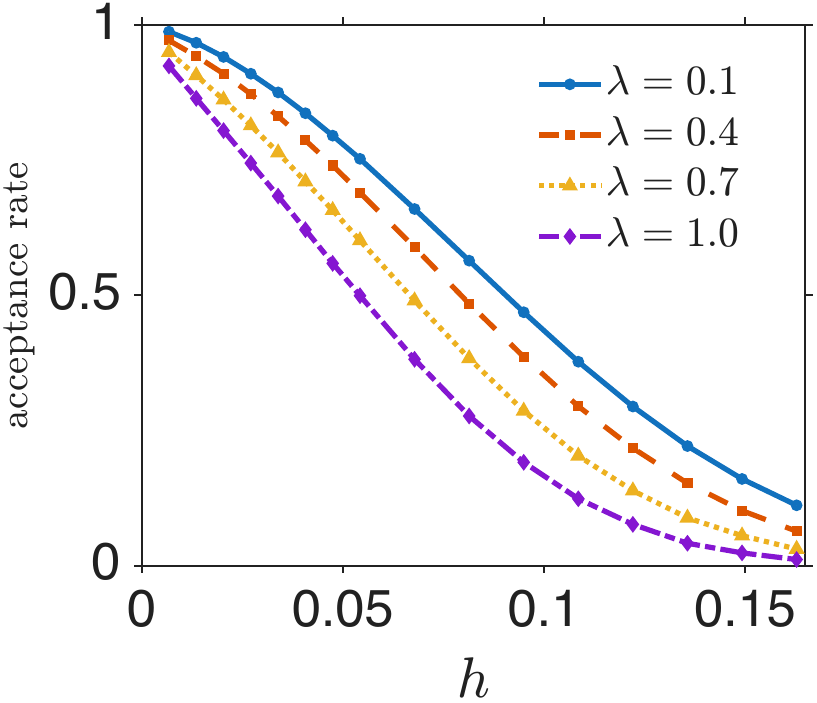}
    \vspace{-5pt}
    \caption{Estimated spectral gap (left) and average acceptance rate (right) of MALA for Bayesian Logistic regression model under different step size $h$ and regularization parameter $\lambda$}
    \label{fig:MALA_LogReg}
\end{figure}

We fix $N=1000$, $d=50$.
We implement the RWM and MALA algorithms to sample from the posterior distribution \eqref{eqn:LogReg}, and estimate the spectral gap following \Cref{sec:Num_BayeLasso}.
We also record the average acceptance rate.
We test under different step size $h$ and regularization parameter $\lambda$.
For each setting, we run $25$ independent trials with chain length $K = 5\times 10^{4}$, preceded by an additional $10^{4}$ burn-in steps, and report averages over trials.
In \Cref{fig:RWM_LogReg,fig:MALA_LogReg}, we can observe similar behaviors as in \Cref{sec:Num_BayeLasso}, which verify our theoretical results in \Cref{sec:SpecGapRWM,sec:SpecGapMALA}.

\section{Conclusion}
\label{sec:Conclusion}

In this work, we established spectral gap bounds of Random--Walk Metropolis and Metropolis--adjusted Langevin algorithms for non-smooth target distributions under convexity-type assumptions.
Our analysis is based on the isoperimetric approach, with a new close coupling condition established for both algorithms under suitable step size constraints.
We also provided numerical experiments to verify our theoretical results.
Future directions include deriving optimal bounds under log-Sobolev or Poincaré inequalities; extending to other Metropolis--Hastings algorithms, including Hamiltonian Monte Carlo; and relaxing the step size constraint for MALA to achieve optimal scaling.

\appendix

\section{Deferred proofs}
\label{app:Proofs}

\subsection{Proof of \texorpdfstring{\Cref{prop:SpecGap_via_Iso}}{Proposition \ref{prop:SpecGap_via_Iso}}}
\label{app:pf_SpecGap_via_Iso}

\begin{proof}
By Cheeger's inequality \Cref{thm:ChIneq}, it suffices to lower bound the Cheeger constant $\Phi(\sfP)$ defined in \eqref{eqn:Cond}.
For any measurable set $S \subseteq \mR^d$, denote the three sets
\[
    S_1 = \left\{ x \in S: \sfP(x,S\matc) \le \epsilon/2 \right\}, \quad S_2 = \left\{ x \in S\matc: \sfP(x,S) \le \epsilon/2 \right\}, \quad S_3 = (S_1\cup S_2)\matc.
\]
By definition, for any $x\in S_1$ and $y \in S_2$,
\[
    \TV \Brac{ \delta_x \sfP , \delta_y \sfP } \ge \sfP(x,S) - \sfP(y,S) \ge 1 - \epsilon/2 - \epsilon/2 = 1 - \epsilon.
\]
By close coupling condition \eqref{eqn:CloseCouple}, this implies that $\norme{x-y} \ge \delta$, and thus $d(S_1,S_2) \ge \delta$.

For any $\theta \in (0,1)$, we discuss by the following two cases: 
\begin{itemize}
    \item Case I: $\pi(S_1) \le \theta \pi(S)$ or $\pi(S_2) \le \theta \pi(S\matc)$.
    \item Case II: $ \pi(S_1) > \theta \pi(S)$ and $ \pi(S_2) > \theta \pi(S\matc)$.
\end{itemize}
(1) Denote for simplicity $\Pi = \pi \otimes \sfP$. If $\pi(S_1) \le \theta \pi(S)$, by definition of $S_1$, we have
\[
    \Pi(S \times S\matc) \ge \Pi((S\setminus S_1) \times S\matc) \ge \frac{\epsilon}{2}  \pi(S\setminus S_1) \ge \frac{\epsilon}{2} (1-\theta) \pi(S).
\]
Since $\sfP$ is reversible with respect to $\pi$, we have $\Pi(S \times S\matc) = \Pi(S\matc \times S)$, and thus the above inequality also holds for $\Pi(S\matc \times S)$.
Likewise, if $\pi(S_2) \le \theta \pi(S\matc)$, we can show that
\[
    \Pi(S\matc \times S) \ge \frac{\epsilon}{2} (1-\theta) \pi(S\matc).
\]
Therefore, in Case I, it always holds that
\[
    \Pi(S \times S\matc) \ge \frac{\epsilon}{2} (1-\theta) \min \{ \pi(S), \pi(S\matc) \}.
\]
(2) When $ \pi(S_1) > \theta \pi(S)$ and $ \pi(S_2) > \theta \pi(S\matc)$, by the three-set isoperimetric inequality,
\[
    \pi(S_3) \ge \Upsilon \Brac{ d(S_1,S_2) } F \Brac{ \min \{ \pi(S_1), \pi(S_2) \} } \ge \Upsilon(\delta)\, F \Brac{ \theta \min \{ \pi(S), \pi(S\matc) \} }.
\]
where the second inequality follows from $d(S_1,S_2) \ge \delta$ and the monotonicity of $\Upsilon$ and $F$.
Now by the definition of $S_1,S_2$, we have
\begin{align*}
    \Pi(S \times S\matc) =~& \frac{1}{2} \Brac{ \Pi(S \times S\matc) + \Pi(S\matc \times S) } \\
    \ge~& \frac{1}{2} \Brac{ \Pi((S\setminus S_1) \times S\matc) + \Pi((S\matc\setminus S_2) \times S) } \\
    \ge~& \frac{\epsilon}{4} \Brac{ \pi(S\setminus S_1) + \pi(S\matc\setminus S_2)  } = \frac{\epsilon}{4} \pi(S_3) \\
    \ge~& \frac{\epsilon}{4} \Upsilon(\delta) F \Brac{ \theta \min \{ \pi(S), \pi(S\matc) \} }.
\end{align*}
Combining the above two cases, we have
\[
    \Pi(S \times S\matc) \ge \min \left\{ \frac{\epsilon}{2} (1-\theta) \min \{ \pi(S), \pi(S\matc) \}, \frac{\epsilon}{4} \Upsilon(\delta) F \Brac{ \theta \min \{ \pi(S), \pi(S\matc) \} } \right\}.
\]
By the symmetry of $S$ and $S\matc$, we can assume without loss of generality that $\pi(S) \le \frac{1}{2}$.
Then
\[
    \frac{\Pi(S \times S\matc)}{\pi(S)} \ge \min \left\{ \frac{\epsilon}{2} (1-\theta), \frac{\epsilon}{4} \Upsilon(\delta) \frac{ F(\theta \pi(S)) }{ \pi(S) } \right\}.
\]
Taking infimum over $S$, we get
\[
    \Phi(\sfP) \ge \frac{\epsilon}{2} \min \left\{ 1-\theta , \Upsilon(\delta) \inf_{t\in(0,\frac{1}{2}]} \frac{ F(\theta t) }{ 2t } \right\}.
\]
The first claim follows by taking supremum over $\theta\in[0,1]$ and applying \eqref{eqn:ChIneq}.

For the second claim, when $F$ is a concave function on $(0,\frac{1}{2}]$, then $\forall 0<x\le y\le \frac{1}{2}$,
\[
    \frac{F(y)-F(x)}{y-x} \le \frac{F(x) - \lim_{t \to 0+} F(t) }{x - 0} \le \frac{F(x)}{x} ~\St~ \frac{F(y)}{y} \le \frac{F(x)}{x}.
\]
Therefore, $\frac{F(t)}{t}$ is non-increasing on $(0,\frac{1}{2}]$, and thus $\inf_{t\in(0,\frac{1}{2}]} \frac{ F(\theta t) }{ 2t } = F(\theta/2)$.
Taking $\theta = \frac{1}{2}$ in \eqref{eqn:SpecGap_via_Iso} yields the desired result.
\end{proof}

\subsection{Proof of \texorpdfstring{\Cref{lem:CloseCoupleRWM}}{Lemma \ref{lem:CloseCoupleRWM}}}
\label{app:pf_CloseCoupleRWM}

The proof of \Cref{lem:CloseCoupleRWM} relies on the following close coupling lemma for Metropolis kernels with symmetric proposals.
Recall that the Metropolis kernel $\sfP$ with proposal kernel $\sfQ$ is given by
\[
    \sfP(x,y) = \alpha(x,y) \sfQ(x,y)
    + r(x) \delta_x(y),
\]
where $\alpha(x,y)$ is the acceptance rate \eqref{eqn:acrate}, and $r(x)$ is the rejection rate
\begin{equation}
\label{eqn:rejectrate}
    r(x) = 1 - \int \alpha(x,y) \sfQ(x,y) \mdd y.
\end{equation}
If the proposal kernel is symmetric, we have the following close coupling lemma.
\begin{lem}[\cite{24-AAP2058}, Lemma 19]
\label{lem:CloseCoupleMH}
Suppose $\sfQ(x,y) = \sfQ(y,x)$. Then for any $x,y$, it holds
\[
    \TV \Brac{ \delta_x \sfP, \delta_y \sfP } \le \TV \Brac{ \delta_x \sfQ, \delta_y \sfQ } + \sup_x r(x).
\]
\end{lem}

Using \Cref{lem:CloseCoupleMH}, we can establish the close coupling condition for Random--Walk Metropolis.

\begin{proof}[Proof of \Cref{lem:CloseCoupleRWM}]
By \Cref{lem:CloseCoupleMH}, since the Random--Walk Metropolis proposal is symmetric, it suffices to control $\TV(\delta_x \sfQ,\delta_y \sfQ)$ and $\sup_x r(x)$. 
First notice $\delta_x \sfQ,\delta_y \sfQ$ are two Gaussians with the same variance, by Pinsker inequality,
\[
    \TV (\delta_x \sfQ,\delta_y \sfQ) \le \sqrt{\frac{1}{2} \KL (\delta_x \sfQ,\delta_y \sfQ)} = \frac{\norme{x-y}}{2h}.
\]
Next we control $\sup_x r(x)$. By \Cref{asm:smth}, we have
\begin{align*}
    &|\log \pi(x+y) - \log \pi(x) - \ip{\nabla f (x)}{y} | \\
    \le~& |f(x+y) - f(x) - \ip{\nabla f (x)}{y} | + |g(x+y) - g(x) | \\
    \le~& \frac{1}{2} \sfM \norme{y}^2 + 2 \sfL \norme{y}.
\end{align*}
Denote $a \wedge b = \min\{a,b\}$ for simplicity.
Let $\psi(t) = 2 \sfL t + \frac{1}{2} \sfM t^2 $.
Then 
\begin{align*}
    1 - r(x) =~& \mE_{\xi\sim \GN(0,\sfI)} \Rectbrac{ \frac{\pi(x+h\xi)}{\pi(x)} \wedge 1 } \\
    =~& \mE_{\xi\sim \GN(0,\sfI)} \Rectbrac{ \mee^{ \log \pi (x+h\xi) - \log \pi(x) } \wedge 1 }  \\
    \ge~& \mE_{\xi\sim \GN(0,\sfI)} \Rectbrac{ \mee^{ h \ip{\nabla f (x)}{\xi} - \psi \Brac{ h \norme{\xi} } } \wedge 1 } \ge \mE_\xi \Rectbrac{ \mee^{-\psi(h\norme{\xi})} \Brac{ \mee^{h\ip{\nabla f (x)}{\xi}} \wedge 1 } }.
\end{align*}
By symmetry of Gaussian distribution, and notice $a\wedge 1 + a^{-1} \wedge 1 \ge 1$, we obtain 
\begin{align*}
    \mE_\xi \Rectbrac{ \mee^{-\psi(h\norme{\xi})} \Brac{ \mee^{h\ip{\nabla f (x)}{\xi}} \wedge 1 } } =~& \frac{1}{2} \mE_\xi \Rectbrac{ \mee^{-\psi(h\norme{\xi})} \Brac{ \mee^{h\ip{\nabla f (x)}{\xi}} \wedge 1 + \mee^{-h\ip{\nabla f (x)}{\xi}} \wedge 1 } } \\
    \ge~& \frac{1}{2} \mE_\xi \Rectbrac{ \mee^{-\psi(h\norme{\xi})} }.
\end{align*}
Therefore, we have 
\begin{align*}
    \sup_x r(x) \le~& 1 - \frac{1}{2} \mE_\xi \Rectbrac{ \mee^{-\psi(h\norme{\xi})} } \\
    \le~& 1 - \frac{1}{2} \mE \Rectbrac{ 1 - \psi(h\norme{\xi}) } \\
    =~& \frac{1}{2} + \frac{1}{2} \mE \Rectbrac{ 2 \sfL h\norme{\xi} + \frac{1}{2} \sfM h^2 \norme{\xi}^2  } \\
    \le~& \frac{1}{2} + \sfL h \sqrt{d} + \frac{1}{4} \sfM h^2 d.
\end{align*}
Take $h\le \frac{1}{\sqrt{d}} \min \left\{ \frac{1}{2\sqrt{\sfM}}, \frac{1}{16\sfL} \right\}$, then $\sup_x r(x) \le \frac{5}{8}$. Finally, when $\norme{x-y}<h/4$, it holds
\[
    \TV \Brac{ \delta_x \sfP, \delta_y \sfP } \le \frac{\norme{x-y}}{2h} + \sup_x r(x) < \frac{3}{4}.
\]
This completes the proof.
\end{proof}

\subsection{Proof of \texorpdfstring{\Cref{lem:ac_MALA}}{Lemma \ref{lem:ac_MALA}}}
\label{app:pf_ac_MALA}

\begin{proof}
Denote MALA proposal as
\[
    x' = x + h v(x) + \sqrt{2h} \xi, \quad \xi \sim \sfN(0,\sfI).
\]
Denote for simplicity $l(x) = \log \pi(x)$.
The acceptance rate is 
\[
    \alpha(x,x') = \min \left\{ 1, \frac{ \pi(x') \sfQ(x',x) }{ \pi(x) \sfQ(x,x') } \right\} = \exp \Brac{ \min \{ 0, a(x,x') \} }, 
\]
where we denote 
\[
    a(x,x') = \log \frac{ \pi(x') \sfQ(x',x) }{ \pi(x) \sfQ(x,x') } = l(x') - l(x) - \frac{1}{4h} \norme{ x - x' - h v(x') }^2 + \frac{1}{4h} \norme{ x' - x - h v(x) }^2. 
\]
Direct calculation shows (recall $x' - x = h v(x) + \sqrt{2h} \xi$)
\[
    a(x,x') = l(x') - l(x) - \ip{ v(x') }{x'-x} - \frac{h}{4} \norme{ v(x') - v(x) }^2 + \frac{\sqrt{2h} }{2} \ip{ \xi }{v(x')-v(x)}. 
\]
Since $l$ is $\sfm$-strongly concave, it holds that
\[
    l(x) \le l(x') + \ip{v(x')}{x-x'} - \frac{m}{2} \norme{ x - x' }^2. 
\]
\[
    \St~ a(x,x') \ge \frac{\sfm}{2} \norme{ x - x' }^2 - \frac{h}{4} \norme{ v(x') - v(x) }^2 + \frac{\sqrt{2h}}{2} \ip{ \xi }{v(x')-v(x)}.   
\]
Recall $v(x) = \nabla f(x) + v_{\rm s}(x)$. By \Cref{asm:smth}, we have
\begin{equation}
\label{eqn:pf_v_lipbound}
    \norme{v(x')-v(x)} \le \norme{ \nabla f(x') - \nabla f(x) } + \norme{v_{\rm s}(x')-v_{\rm s}(x)} \le \sfM \norme{x'-x} + 2 \sfL.
\end{equation}
So that if $h \le \sfm/(2\sfM^2)$, we have 
\[
    \frac{h}{4} \norme{ v(x') - v(x) }^2 \le \frac{h}{4} \Brac{ 2 \sfM^2 \norme{x'-x}^2 + 8 \sfL^2 } \le \frac{\sfm}{4} \norme{x'-x}^2 + 2 h \sfL^2.
\]
\begin{equation}
\label{eqn:pf_a_lowerbound}
    \St~ a(x,x') \ge \frac{\sfm}{4} \norme{x'-x}^2 - 2 h \sfL^2 + \frac{\sqrt{2h}}{2} \ip{ \xi }{v(x')-v(x)}. 
\end{equation}
Next we will discuss by two cases: $\norme{v(x)} \le R$ or $\norme{v(x)}> R$, where
\begin{equation}
\label{eqn:pf_R_def}
    R := \frac{2(2\sfm+3\sfM)}{\sfm} \frac{\sqrt{d}+\sqrt{6}}{\sqrt{2h}}.
\end{equation}

\noindent
{\bf Case 1: $\norme{v(x)} \le R$.} 
Recall $ x' - x = h v(x) + \sqrt{2h} \xi $, so that by \eqref{eqn:pf_v_lipbound},
\[
    \norme{v(x')-v(x)} \le \sfM \norme{x'-x} + 2 \sfL \le \sfM h \norme{v(x)} + \sfM \sqrt{2h} \norme{\xi} + 2 \sfL.
\]
Therefore, by \eqref{eqn:pf_a_lowerbound}, we have
\begin{align*}
    a(x,x') \ge~& - 2 h \sfL^2 + \frac{\sqrt{2h}}{2} \ip{ \xi }{v(x')-v(x)} \\
    \ge~& - 2 h \sfL^2 - \frac{\sqrt{2h}}{2} \Brac{ \sfM h \norme{v(x)} + \sfM \sqrt{2h} \norme{\xi} + 2 \sfL } \norme{ \xi }. 
\end{align*}
Since $\exp(x)$ is $1$-Lipschitz on $(-\infty,0]$, we have (denote $r(x)$ as the rejection rate as in \eqref{eqn:rejectrate})
\begin{align*}
    r(x) := ~& 1 - \mE_{\xi\sim \GN(0,\sfI)} \Rectbrac{ \exp(\min\{0,a(x,x')\}) } \\
    \le~& \mE_{\xi\sim \GN(0,\sfI)} \Rectbrac{ 1 - \exp \Brac{ - 2 h \sfL^2 - \frac{\sqrt{2h}}{2} \Brac{ \sfM h \norme{v(x)} + \sfM \sqrt{2h} \norme{\xi} + 2 \sfL } \norme{ \xi } } } \\
    \le~& \mE_{\xi\sim \GN(0,\sfI)} \Rectbrac{ 2 h \sfL^2 + \frac{\sqrt{2h}}{2} \Brac{ \sfM h \norme{v(x)} \norme{ \xi } + \sfM \sqrt{2h} \norme{\xi}^2 + 2 \sfL \norme{ \xi } }  } \\
    \le~& 2 h \sfL^2 + \frac{\sqrt{2h}}{2} \Brac{ \sfM h R \sqrt{d} + \sfM \sqrt{2h} d + 2 \sfL \sqrt{d} }.  
\end{align*}
Denote $\kappa = \sfM/\sfm \ge 1$.
When $h \le \frac{1}{200} \cdot \min \left\{ \frac{1}{d\kappa \sfM} , \frac{1}{d\sfL^2} \right\}$,
it holds that (recall the definition of $R$ in \eqref{eqn:pf_R_def})
\begin{align*}
    r(x) \le~& 2 h \sfL^2 + \frac{1}{2} \sfM h R \sqrt{2hd} + \sfM h d + \sfL \sqrt{2hd} \\
    \le~& \frac{1}{100d} + \frac{1}{200} \cdot \Brac{ \frac{2}{\kappa} + 3 } \Brac{ 1 + \sqrt{\frac{6}{d}} } + \frac{1}{200\kappa} + \sqrt{\frac{1}{100}} \\ 
    \le~& 0.01 + 0.005 \cdot \Brac{ 2 + 3 } \Brac{ 1 + \sqrt{6} } + 0.005 + \sqrt{0.01} < 0.35. 
\end{align*}

\noindent
{\bf Case 2: $\norme{v(x)} > R$.} 
Recall $ x' - x = h v(x) + \sqrt{2h} \xi $, so that by \eqref{eqn:pf_a_lowerbound} and \eqref{eqn:pf_v_lipbound}
\begin{equation}
\label{eqn:pf_a_case2}
\begin{split}
    a(x,&\,x') \ge \frac{\sfm}{4} \norme{ h v(x) + \sqrt{2h} \xi }^2 - 2 h \sfL^2 - \frac{\sqrt{2h}}{2} \norme{ \xi } \norme{v(x')-v(x)} \\
    =~& \frac{\sfm h^2}{4} \norme{v(x)}^2 + \frac{\sfm h}{2} \norme{\xi}^2 + \frac{\sfm h\sqrt{2h}}{2} \ip{v(x)}{\xi} - 2 h \sfL^2 - \frac{\sqrt{2h}}{2} \norme{ \xi } \norme{v(x')-v(x)} \\
    \ge~&\frac{\sfm h^2}{4} \norme{v(x)}^2 + \frac{\sfm h}{2} \norme{\xi}^2 - \frac{\sfm h\sqrt{2h}}{2} \norme{v(x)} \norme{\xi} - 2 h \sfL^2 - \frac{\sqrt{2h}}{2} \Brac{ \sfM h \norme{v(x)} + \sfM \sqrt{2h} \norme{\xi} + 2 \sfL } \norme{\xi} \\
    =~& \frac{\sfm h^2}{4} \norme{v(x)}^2 - \frac{(\sfm+\sfM) h\sqrt{2h}}{2} \norme{v(x)} \norme{\xi} - \Brac{ \sfM - \frac{\sfm}{2}  } h \norme{\xi}^2 - \sqrt{2h} \sfL \norme{\xi} - 2 h \sfL^2.
\end{split}
\end{equation}
Consider a high probability event $\{ \norme{\xi} \le \sqrt{d} + \sqrt{6} \}$.
Since $\norme{v(x)} > R$, we have 
\[
    \norme{v(x)} > \frac{2(2\sfm+3\sfM)}{\sfm} \frac{\sqrt{d}+\sqrt{6}}{\sqrt{2h}} \ge \frac{2(2\sfm+3\sfM)}{\sfm} \frac{\norme{\xi}}{\sqrt{2h}}. 
\]
Therefore, using the above inequality, when $\norme{\xi} \le \sqrt{d} + \sqrt{6}$, it holds that
\begin{align*}
    & \frac{(\sfm+\sfM) h\sqrt{2h}}{2} \norme{v(x)} \norme{\xi} + \Brac{ \sfM - \frac{\sfm}{2}  } h \norme{\xi}^2 \\
    \le~& \frac{(\sfm+\sfM) h\sqrt{2h}}{2} \norme{v(x)} \frac{\sfm\sqrt{2h}}{2(2\sfm+3\sfM)} \norme{v(x)} + \Brac{ \sfM - \frac{\sfm}{2}  } h \cdot \Brac{ \frac{\sfm\sqrt{2h}}{2(2\sfm+3\sfM)} \norme{v(x)} }^2 \\
    =~& \frac{ (\sfm+\sfM) \sfm h^2}{2(2\sfm+3\sfM)} \norme{v(x)}^2 + \frac{ \Brac{ 2 \sfM - \sfm } \sfm^2 h^2}{4(2\sfm+3\sfM)^2} \norme{v(x)}^2 = \Brac{ \frac{ 2(\sfm+\sfM) }{2\sfm+3\sfM} + \frac{ \Brac{ 2 \sfM - \sfm } \sfm }{(2\sfm+3\sfM)^2} } \frac{\sfm h^2}{4} \norme{v(x)}^2.
\end{align*}
Now notice
\begin{align*}
    \frac{ 2(\sfm+\sfM) }{2\sfm+3\sfM} + \frac{ \Brac{ 2 \sfM - \sfm } \sfm }{(2\sfm+3\sfM)^2} =~& 1 - \frac{ \sfM }{2\sfm+3\sfM} + \frac{ \Brac{ 2 \sfM - \sfm } \sfm }{(2\sfm+3\sfM)^2}  \\
    =~& 1 - \frac{ \sfM (2\sfm+3\sfM) -  \Brac{ 2 \sfM - \sfm } \sfm }{(2\sfm+3\sfM)^2} = 1 - \frac{\sfm^2+3\sfM^2}{(2\sfm+3\sfM)^2} \le 1.
\end{align*}
Thus, we have
\[
    \frac{(\sfm+\sfM) h\sqrt{2h}}{2} \norme{v(x)} \norme{\xi} + \Brac{ \sfM - \frac{\sfm}{2}  } h \norme{\xi}^2 \le \frac{\sfm h^2}{4} \norme{v(x)}^2.
\]
Plugging the above inequality into \eqref{eqn:pf_a_case2}, we get
\[
    a(x,x') \ge - \sqrt{2h} \sfL \norme{\xi} - 2 h \sfL^2.
\]
So that when $\norme{\xi} \le \sqrt{d} + \sqrt{6} $, it holds that $a(x,x') \ge - \sqrt{2h} \sfL \norme{\xi} - 2 h \sfL^2$, which implies
\begin{align*}
    &\mE_{\xi\sim \GN(0,\sfI)} \Rectbrac{ \exp \Brac{ \min \{ 0, a(x,x') \} } } \\
    \ge~& \mE_{\xi\sim \GN(0,\sfI)} \Rectbrac{ \exp \Brac{ - \sqrt{2h} \sfL \norme{\xi} - 2 h \sfL^2 } \ind_{\norme{\xi} \le \sqrt{d} + \sqrt{6}} } \\
    \ge~& \mP \Brac{ \norme{\xi} \le \sqrt{d} + \sqrt{6} } \cdot \exp \Brac{ - \sqrt{2h} \sfL ( \sqrt{d} + \sqrt{6} ) - 2 h \sfL^2 }.
\end{align*}
By the concentration inequality for $\GN(0,\sfI)$,
\[
    \mP \Rectbrac{ \norme{\xi} \ge \sqrt{d} + \sqrt{6} } \le \mee^{-3}.
\]
When $h \le 1/(200 d\sfL^2)$, it holds that
\begin{align*}
    r(x) =~& 1 - \mE_{\xi\sim \GN(0,\sfI)} \Rectbrac{ \exp \Brac{ \min \{ 0, a(x,x') \} } } \\
    \le~& 1 - (1-\mee^{-3}) \cdot \exp \Brac{ - 0.1 \cdot ( 1 + \sqrt{6/d} ) - 0.01/d} \\
    \le~& 1 - (1-\mee^{-3}) \cdot \exp \Brac{ - 0.1 \cdot ( 1 + \sqrt{6} ) - 0.01 } < 0.35. 
\end{align*}
Combining the two cases, we conclude that when $h \le \frac{1}{200} \cdot \min \{ \frac{1}{d\kappa \sfM} , \frac{1}{d\sfL^2} \} $, it holds that
\[
    \inf_x \mE_{y \sim \sfQ(x,y)} [\alpha(x,y)] = 1 - \sup_x r(x) > 1 - 0.35 = \frac{13}{20}.
\]
This completes the proof.
\end{proof}

\subsection{Proof of \texorpdfstring{\Cref{lem:CloseCoupleMALA}}{Lemma \ref{lem:CloseCoupleMALA}}}
\label{app:pf_CloseCoupleMALA}

\begin{proof}
We first prove an analogue of \Cref{lem:CloseCoupleMH}.
Let $r(x)$ be the rejection rate defined in \eqref{eqn:rejectrate}.
Then by definition of MALA proposal, we have for any $x$,
\[
    \TV(\delta_x \sfP, \delta_x \sfQ) = \frac{1}{2} r(x) + \frac{1}{2} \int ( \sfQ(x,z) - \alpha(x,z) \sfQ(x,z) ) \mdd z = r(x). 
\]
Therefore, by triangle inequality,
\begin{align*}
    \TV \Brac{ \delta_x \sfP, \delta_y \sfP } \le~& \TV \Brac{ \delta_x \sfQ, \delta_y \sfQ } + \TV(\delta_x \sfP, \delta_x \sfQ) + \TV(\delta_y \sfP, \delta_y \sfQ) \\
    =~& \TV \Brac{ \delta_x \sfQ, \delta_y \sfQ } + r(x) + r(y) \le \TV \Brac{ \delta_x \sfQ, \delta_y \sfQ } + 2 \sup_x r(x).
\end{align*}
Thus, it suffices to control $\TV(\delta_x \sfQ,\delta_y \sfQ)$ and $\sup_x r(x)$. 
Since $\delta_x \sfQ,\delta_y \sfQ$ are two Gaussians with the same variance, by Pinsker inequality,
\[
    \TV (\delta_x \sfQ,\delta_y \sfQ) \le \sqrt{\frac{1}{2} \KL (\delta_x \sfQ,\delta_y \sfQ)} = \frac{\norme{x+hv(x)-y-hv(y)}}{2\sqrt{2h}}.
\]
Notice 
\[
    \norme{x+hv(x)-y-hv(y)} \le \norme{x+h\nabla f (x)-y-h\nabla f (y)} + h \norme{v_{\rm s}(x)-v_{\rm s}(y)}. 
\]
For the first term, note $0 \preceq - \nabla^2 f \preceq \sfM I$ by \Cref{asm:MALA}. Since $h < \sfM^{-1}$, it holds that
\[
    0 \preceq \mean{H} := \int_0^1 \Brac{ \sfI + h \nabla^2 f (tx+(1-t)y) } \mdd t \preceq  \sfI.
\]
So that by Taylor's theorem, we have
\[
    \norme{x+h\nabla f (x)-y-h\nabla f (y)} = \norme{\mean{H} (x-y)} \le \norme{x-y}. 
\]
For the second term, by \Cref{asm:smth}, $\norme{v_{\rm s}(x)-v_{\rm s}(y)} \le 2 \sfL$.
Combined, we have
\[
    \norme{x+hv(x)-y-hv(y)} \le \norme{x-y} + 2 h \sfL. 
\]
When $h\le \min \{ \frac{1}{\sfM} , \frac{1}{200\sfL^2d} \}$, and $\norme{x-y} < \frac{\sqrt{2h}}{10}$, it holds that
\[
    \TV (\delta_x \sfQ,\delta_y \sfQ) \le \frac{\norme{x-y} + 2 h \sfL}{2\sqrt{2h}} \le \frac{1}{20} + \frac{1}{20\sqrt{d}} \le \frac{1}{10}.  
\]
Finally, by \Cref{lem:ac_MALA},
\[
    \sup_x r(x) = 1 - \inf_x \E_{y \sim \sfQ(x,y)} [\alpha(x,y)] < 1 - \frac{13}{20} = \frac{7}{20}.
\]
\[
    \St\quad \TV(\delta_x \sfP,\delta_y \sfP) \le
    \TV(\delta_x \sfQ,\delta_y \sfQ) + 2 \sup_x r(x) < \frac{1}{10} + \frac{7}{10} = \frac{4}{5}.
\]
This completes the proof. 
\end{proof}

\section*{Acknowledgements}

The work of SL and XT is partially supported by the National University of Singapore and Singapore MOE grant Tier-1-A-8002956-00-00.

\bibliographystyle{apalike}
\bibliography{BibLib}

\end{document}